\documentclass[11pt, leqno, a4paper]{amsart}
\usepackage{amsmath,amstext,amsthm,amsfonts}
\usepackage{amssymb,latexsym}
\usepackage[cp1252]{inputenc}
\usepackage[T1]{fontenc}
\usepackage{lmodern}
\usepackage{hyperref, color}
\usepackage[english]{babel}
\usepackage{mathrsfs}

\newtheorem{theorem}{\textbf{Theorem}}[section]
\newtheorem{proposition}[theorem]{\textbf{Proposition}}
\newtheorem{prop}[theorem]{\textbf{Proposition}}
\newtheorem{thm}[theorem]{\textbf{Theorem}}

\newtheorem{lem}[theorem]{\textbf{Lemma}}

\theoremstyle{definition}
\theoremstyle{remark}

\newcommand{\field}[1]{\mathbb{#1}}
\newcommand{\real}{\field{R}} %%%%%% %%%%% reais

\newcommand{\al} {\alpha}       
\newcommand{\be} {\beta}        
\newcommand{\ga} {\gamma}    
       \newcommand{\De}{\Delta}

\newcommand{\te} {\theta}

\newcommand{\la} {\lambda}

       \newcommand{\Om}{\Omega}

\newcommand{\ie}{\emph{i.e. }}

\newcommand{\HH}{\mathbb{H}}%%%%%%%%
\newcommand{\R}{\mathbb{R}}
\newcommand{\N}{\mathbb{N}}

\newcommand{\cF}{\mathcal{F}}
\newcommand{\cE}{\mathcal{E}}
\newcommand{\cK}{\mathcal{K}}
\newcommand{\cC}{\mathcal{C}}
\newcommand{\cD}{\mathcal{D}}

\newcommand{\Mh}{\widehat{M}}

\newcommand{\rich}{\widehat{\mathrm{Ric}}}

\newcommand{\gh}{\hat{g}}

\newcommand{\pf}{\par{\noindent\textbf{Proof.~}}}
\newcommand{\nil}{\mathrm{Nil}(3)}       %   Nil3, using mathrm
\newcommand{\niln}{\mathrm{Nil}(2n+1)}
\newcommand{\rh}{\widehat{\mathrm{Ric}}} %   Ricci, using mathrm

\newcommand{\dd}{\mathrm{d}}             %   special differential d, with mathrm

\begin{document}

%%%%%%%%%%%%%%%%%%%%%%%%%%%%%%%%%%%%%%%%%%%%%%%%%%%%%%%%%%

\title[Rotational Catenoids in the Heisenberg Groups]{Stability
Properties of Rotational Catenoids in the Heisenberg Groups}

\author[P. B\'erard, M.  Cavalcante]{Pierre B\'erard,  Marcos P. Cavalcante}

%\date{\today}
\date{March 13, 2013. Final version to appear in Matem\'{a}tica Contempor\^{a}nea \\(Sociedade
Brasileira de Matem\'{a}tica)} \maketitle

\thispagestyle{empty}

%%%%%%%%%%%%%%%%%%%%%%%%%%%%%%%%%%%%%

\begin{abstract}
\noindent {In this paper, we determine the maximally stable,
rotationally invariant domains on the catenoids $\cC_a$ (minimal
surfaces invariant by rotations) in the Heisenberg group with a
left-invariant metric. We show that these catenoids have Morse index
at least $3$ and we bound the index from above in terms of the
parameter $a$. We also show that the index of $\cC_a$ tends to
infinity with $a$. Finally, we study the rotationally symmetric
stable domains on the higher dimensional catenoids}.
\end{abstract}\bigskip

\textbf{MSC}(2010): 53C42, 58C40.\bigskip

\textbf{Keywords}: Minimal Surface, Heisenberg Group, Killing Field,
Index.

\bigskip

%%%%%%%%%%%%%%%%%%%%%%%%%%%%%%%%%%%%%
\section{Introduction}\label{s-in}

Minimal surfaces in the Heisenberg group equiped with a
left-invariant metric have been studied by several authors, see
\cite{Fi96, FMP99, Dan11, FM09} and the references therein.
\emph{Catenoids} in the Heisenberg group $\nil$ are complete minimal
surfaces which are invariant under a one-parame\-ter subgroup of
rotations with axis the center of the group. They come in a
one-parameter family $\{\mathcal{C}_a, a > 0\}$ of complete minimal
surfaces and were first described in \cite{Fi96} and \cite{FMP99}
where the authors provide the classification of constant mean
curvature surfaces in the Heisenberg group, invariant under certain
subgroups of isometries (the parameter $a$ is the neck size of the
catenoid, see \eqref{cp}).
\medskip

In this paper, we study the stability properties of the catenoids
$\{\cC_a, a> 0\}$. More precisely, we determine the rotationally
invariant stable domains of the catenoids in $\niln$, $n \ge 1$,
with a different behaviour (Lindeloef's property) when $n=1$ and
when $n\ge 2$. We also study the Morse index of the catenoids in
$\nil$. As in \cite{BS10}, the proofs rely in part on a detailed
analysis of the Jacobi fields induced from the Killing fields of the
ambient Heisenberg space and from the variation of the parameter
$a$.\medskip

The paper is organized as follows. In Section \ref{s-pc}, we give
some preliminary results. We first recall the basic geometry of the
Heisenberg group $\nil$ equiped with a left-invariant metric $\gh$
(see \cite{FMP99} for more details). In order to keep our paper
self-contained, we derive the differential equation satisfied by the
generating curves of the catenoids, using a flux formula. In Section
\ref{s-sd}, we describe the stable rotationally invariant domains on
$\{\mathcal{C}_a\}$ (Theorem~\ref{t1}). The proof uses Jacobi
fields. We also give some information on the Gauss map of the
catenoids $\{\mathcal{C}_a\}$. In Section \ref{s-ind},
Theorem~\ref{t2}, we prove that the catenoids $\cC_a, a > 0$ have
Morse index at least $3$. We bound the index from above in terms of
$a$, and we also show that its goes to infinity with $a$. The proof
uses Jacobi fields, Fourier analysis and an adapted perturbation of
the original parametrization of the catenoids. Finally, in
Section~\ref{s-hd}, we study the maximally stable, rotationally
invariant domains on the higher dimensional catenoids
(Theorem~\ref{t3}).\medskip

In the sequel our functions will often depend on the parameter $a$.
We will occasionally omit $a$ to keep the notations simpler. In this
paper, we only consider left-invariant Riemannian metrics on the
Heisenberg groups. \bigskip

The first author was partially supported by the cooperation
programme Math-AmSud. The second author would like to thank Institut
Fourier (Grenoble) for their hospitality during the preparation of
this paper. He gratefully acknowledges CAPES and CNPq for their
financial support.

%%%%%%%%%%%%%%%%%%%%%%%%%%%%%%%%%%%%%%%%%%%%%%%%%%%%%%%%%%%%%%%%%%

\section{Preliminaries}\label{s-pc}

\subsection{The $3$-dimensional Heisenberg manifold}

Let $\nil$ denote the $3$-dimensional \emph{Heisenberg group}. This
is a two-step nilpotent Lie group which can be seen as the subgroup
of $3 \times 3$ matrices given by
\begin{displaymath}
\nil = \left\{ \left( \begin{array}{ccc}
1 & x & z \\
0 & 1 & y \\
0 & 0 & 1
\end{array} \right);\, (x,y,z)\in\real^3 \right\} \subset GL(3,\real).
\end{displaymath}

We denote the corresponding Lie algebra by
$$
\mathcal{L}(\nil) = \left\{ \left( \begin{array}{ccc}
0 & x & z \\
0 & 0 & y \\
0 & 0 & 0
\end{array} \right);\, (x,y,z)\in\real^3 \right\}.
$$

Using the exponential map, $\exp : \mathcal{L}(\nil) \to \nil$, and
the Campbell-Hausdorff formula,
$$
\exp\big(A\big) \exp\big(B\big) = \exp
\big(A+B+\frac{1}{2}[A,B]\big), ~~\forall A, B \in
\mathcal{L}(\nil),
$$
we can view $\nil$ as $\real^3$ equiped with the group structure
$\star$ given by
\begin{equation*}\label{E-star}
(x,y,z) \star (x',y',z') = \Big( x+x', y+y', z+z'+ \frac{1}{2}(xy' -
x'y) \Big),
\end{equation*}
with neutral element $0 = (0,0,0)$ and inverse $\check{p}$ of $p =
(a,b,c)$ given by $\check{p} = (-a,-b,-c)$. The left-multiplication
by $p$ in $\nil$, $L_p : q \mapsto p \star q$, has tangent map
\begin{equation}\label{E-tg}
T_qL_p = \begin{pmatrix}
  1 & 0 & 0 \\
  0 & 1 & 0 \\
  - \frac{1}{2}b & \frac{1}{2}a & 1 \\
\end{pmatrix}
\end{equation}

in the canonical coordinates $\{x,y,z\}$ of $\real^3$ (they are
often referred to as \emph{exponential coordinates}). Let
$\{\partial_x,
\partial_y, \partial_z\}$ denote the cano\-nical vector fields in $\real^3$. It
follows from the expression (\ref{E-tg}) that the vector fields
\begin{equation}\label{E-liv}
\left\{%
\begin{array}{lllll}
X(x,y,z) & = & T_0 L_{(x,y,z)}(\partial_x) &=& \partial_x -
\frac{y}{2}\, \partial_z, \\[4pt]
Y(x,y,z) & = & T_0 L_{(x,y,z)}(\partial_y) &=& \partial_y +
\frac{x}{2}\, \partial_z, \\[4pt]
Z(x,y,z) & = & T_0 L_{(x,y,z)}(\partial_z) &=& \partial_z, \\
\end{array}
\right.
\end{equation}
form a basis of left-invariant vector fields in $\nil$. \medskip

\textbf{The metric $\gh$ on $\nil$}. From now on, we fix the
left-invariant metric $\gh$ on $\nil$ to be such that the family
$\{X,Y,Z\}$ is an orthonormal frame. In the coordinates $\{x,y,z\}$,
this metric is given by
$$
\gh = \dd x^2 + \dd y^2 + \big( \dd z + \frac{1}{2}(y \, \dd x - x
\, \dd y)\big)^2.
$$

The following properties are well-known and can be found for example
in \cite{FMP99}, Section 1. Equiped with the left-invariant metric
$\gh$, the Heisenberg group $\nil$ is a homogeneous Riemannian
manifold whose group of isometries has dimension $4$. A basis of
Killing vector fields on $(\nil, \gh)$ is given by
\begin{equation*}\label{killing}
\left\{%
\begin{array}{cll}
\xi &= & X + y Z,\\
\eta &= & Y -  x  Z,\\
\zeta &= & Z,\\
\rho &= & yX - xY + \frac{1}{2}(x^2 + y^2)Z.
\end{array}
\right.%
\end{equation*}

The first three vector fields $\xi, \eta$ and $ \zeta$ correspond to
the one-parameter subgroups of isometries generated by
right-invariant vector fields in $\nil$, while the vector field
$\rho$ corresponds to the one-parameter subgroup of isometries
defined by
\begin{equation}\label{E-rot}
\psi_{\theta}\big( (x,y,z) \big) = \big( x \cos \theta - y \sin
\theta, x \sin \theta + y \cos \theta,z \big), (x,y,z) \in \R^3,
\theta \in \R,
\end{equation}
in the representation $(\real^3, \star)$ of $\nil$. We call them
\emph{rotations around the $z$-axis}. Notice that the $z$-axis is
precisely the center of $\nil$.

%%PRIVATE

\subsection{Surfaces of revolution in $\nil$}

We say that a surface $M$ in $\nil$ is a \emph{surface of
revolution} if $M$ is invariant under the action of the
one-parameter subgroup $\{\psi_{\theta}, \theta \in \R\}$ given by
(\ref{E-rot}). We will consider surfaces of revolution whose
generating curves are graphs $t \to \big( f(t),t \big)$ above the
$z$-axis in the $2$-plane $\{x,z\}$,  where $f$ is a positive
function, and where $t$ varies in some interval $I \subset \R$. They
are given by a map
\begin{equation}\label{param}
\mathcal{F}(t,\theta) = (f(t)\cos \te, f(t)\sin\te, t),
\end{equation}
for $t \in I \subset \real$ and $\te \in [0, 2 \pi]$. \medskip

\emph{Catenoids}, \ie minimal surfaces of revolution, in $\nil$ are
described in \cite{Fi96, FMP99}, using the methods of equivariant
differential geometry. They come in a one-parameter family of
complete minimal surfaces, $\{\cC_a, a>0\}$. For the sake of
completeness and for later purposes, we will derive the differential
equations satisfied by the generating curve of a catenoid using a
\emph{flux formula} which we now state.

\begin{proposition}\label{P-flux}
Let $(M^n,g) \looparrowright (\Mh^{n+1}, \gh)$ be an isometric
immersion with Riemannian measure $\mu_g$ and normalized mean
curvature vector $\vec{H}$. Let $\Omega$ be a relatively compact
smooth domain in $M$. Let $ \nu_{int}$ denote the unit normal to
$\partial \Omega$ in $M$, pointing inwards, and $\sigma_{g}$ the Riemannian
measure on $\partial \Omega$ induced by $g$. Then, for any Killing
vector field $\mathcal{K}$ on $\Mh^{n+1}$, we have
\begin{equation*}\label{E-flux}
\int_{\partial \Omega} \gh(\mathcal{K}, \nu_{int})\, \dd \sigma_{g}=
- n \int_\Omega \gh(\mathcal{K}, \vec{H}) \, \dd \mu_g.
\end{equation*}
\end{proposition}\medskip

\pf Let $\kappa$ be the restriction to $M$ of the $1$-form dual to
$\mathcal{K}$, i.e. $\kappa =\gh (\mathcal{K}, .)|_M$. Recall that
$\mathcal{K}$ is a Killing field if and only if, for any vector
field $X$ on $\Mh$,  $\gh (\widehat{D}_X{\mathcal{K}},X)=0$
(\cite{KN63}, Proposition~3.2, p. 237). A straightforward
computation shows that the divergence $\delta_g \kappa$ of the
1-form $\kappa$, for the induced metric $g$ on $M$, is given by
$$
\delta_g \kappa= -n\gh(\mathcal{K}, \vec{H}).
$$

The proposition follows from the divergence theorem. \qed
\medskip

%%PRIVATE

Let $M = \cF (I \times [0,2\pi])$ be a minimal surface of revolution
in $\nil$, given by an immersion $\mathcal{F}(t,\te)$ as in
(\ref{param}), with $t \in I \subset \real, \, \te \in[0,2\pi].$ We
can make a coherent choice of a unit vector field $\nu$ tangent to
$M$ and orthogonal to the circles $C_t = \cF(\{t\} \times [0,2\pi])$
in such a way that Proposition~\ref{P-flux} gives
\begin{equation}\label{flux-c}
\int_{C_t}\gh(\mathcal{K},\nu)\, \dd \sigma_{C_t}=
\int_{C_{t_0}}\gh(\mathcal{K},\nu)\, \dd \sigma_{C_{t_0}},
\end{equation}
for all $t_0, t \in I$ and for any Killing vector field
$\mathcal{K}$ in $\nil$.\medskip

\begin{proposition}\label{P-cat}
The generating curve of a minimal surface of revolution of
the form (\ref{param}) in $\nil$ satisfies the first order differential
equation
\begin{equation}\label{E-cat1}
f \big( 4 + f^2  f_t^2 + 4 f_t^2 \big)^{-1/2} = C \text{~~(a
constant)},
\end{equation}
and the second order differential equation
\begin{equation}\label{E-cat2}
f (4+f^2) \, f_{tt} = 4 \, (1+f_t^2),
\end{equation}
where $f_t$ and $f_{tt}$ denote respectively the first and second
derivatives of the function $f$ with respect to the variable $t$.
\end{proposition}

\pf According to \cite{FMP99} Theorem 3, we already know that
minimal surfaces of revolution do exist in $\nil$. Equation
(\ref{E-cat1}) is established by applying Proposition~\ref{P-flux}
with the Killing field $\mathcal{K} = Z$. The constant $C$ can then
be interpreted in terms of a flux. The vectors $\cF_t$ and
$\cF_{\theta}$ are tangent to the surface. Using (\ref{E-liv}), they
can be expressed in the orthonormal frame $\{X,Y,Z\}$ at
$\cF(t,\theta)$ as
\begin{equation}\label{tangents}
\left\{%
\begin{array}{cll}
\cF_t &=& f_t \cos \theta \, X + f_t \sin \theta \, Y + Z, \\[4pt]
\cF_{\theta} &=& - f \sin \theta \, X + f \cos \theta \, Y -
\frac{1}{2} f^2 \, Z.
\end{array}%
\right.
\end{equation}

The Riemannian measure $\sigma_{C_t}$ is given by
$$
\dd \sigma_{C_t} = \sqrt{\gh(\mathcal{F}_{\theta},
\mathcal{F}_{\theta})}\, \dd \theta = f \sqrt{1+\frac{1}{4}f^2} \,
\dd \theta.
$$
Up to sign, the vector $\nu$ is characterized by the facts that it
is unitary, tangent to the surface -- hence a linear combination of
$\mathcal{F}_t$ and $\mathcal{F}_{\theta}$ -- and orthogonal to
$\mathcal{F}_{\theta}$. Consider the vector $n = \cF_t + \alpha
\cF_{\theta}$ with $\alpha$ such that $\gh(n,\cF_{\theta}) = 0$.
Choose $\nu = \gh(n,n)^{-1/2} n$. The expression $\gh(Z,\nu)$ which
appears in (\ref{flux-c}) when we choose $\cK=Z$, is the
$Z$-component of $\nu$. A straightforward computation gives that
$\alpha = 2 (4+f^2)^{-1}$, $\gh(n,Z) = 4 (4+f^2)^{-1}$ and $\gh
(n,n) = f_t^2 + 4(4+f^2)^{-1}$. It follows that
$$\gh(Z,\nu) = 4 (4+f^2)^{-1} \big( f_t^2 + \frac{4}{4+f^2}\big)^{-1/2}.$$

Using (\ref{flux-c}),we obtain that the quantity (a flux)
$$
f(t) \big[ 4 + f^2(t)\, f_t^2(t) + 4 f_t^2(t) \big]^{-1/2}
$$
is independent of $t$. Equation~(\ref{E-cat1}) follows. Taking the
derivative of (\ref{E-cat1}) and using the fact that $f_t \not
\equiv 0$ (see \cite{FMP99}), we obtain Equation~(\ref{E-cat2}).
\qed \medskip

\textbf{Remark}. The above equations can also be derived directly
from \cite{FMP99} (using the computations in the proof of their
Theorem~3) or by minimizing the area of a rotational domain, in the
spirit of the calculus of variations.

\subsection{Qualitative analysis of Equation (\ref{E-cat2})}\label{ss-qa}
Given $a>0$, consider the \emph{Cauchy problem},
\begin{equation}\label{cp}
\left\{%
\begin{array}{rcl}
f (f^2+4) f_{tt} &=& 4(1+f_t^2),\\
f(0)&=&a, \\
f_t(0)&=&0,
\end{array}%
\right.
\end{equation}
where the subscript $t$ means that we take the derivative with
respect to $t$. Recall that this differential equation admits a
first integral and, more precisely, that
\begin{equation}\label{1st}
\frac{(f^2+4)(1+f_t^2)}{f^2}= \frac{a^2+4}{a^2}.
\end{equation}
%
%%PRIVATE
A simple analysis shows that (\ref{cp}) admits a maximal solution
$f(a,t)$ which is an even function of $t$ on some interval $(-A_a,
A_a)$. Furthermore, the function
$$
f(a,\cdot):[0,A_a)\to [a,\infty)
$$
is an increasing function and we can introduce its inverse function
$$
\phi(a,\cdot):[a,\infty)\to[0,A_a).
$$

Using (\ref{1st}), we infer that $\phi$ is given by the integral
\begin{equation}\label{E-phi}
\phi(a,\tau)= \frac{a}{2}\int_1^{\tau/a}
\sqrt{\frac{a^2v^2+4}{v^2-1}}\,\dd v.
\end{equation}
It follows that
\begin{equation*}\label{assint}
\phi(a,\tau) \thicksim \frac{a}{2}\tau, \text{ ~when~ }\tau \to
\infty.
\end{equation*}

Finally, we conclude that the Cauchy problem (\ref{cp}) admits a
global solution $f(a, \cdot): \real \to [a,\infty)$ which satisfies
\begin{equation*}\label{gs}
\left\{%
\begin{array}{l}
f(a,t)=f(a,-t),\\[5pt]
f(a,t) \thicksim \frac{2}{a} |t|, \quad \text{and}\\[5pt]
f_t(a,t) \thicksim \frac{2}{a} \mathrm{sgn}(t), \quad \text{when~ }
|t| \to \infty.
\end{array}%
\right.
\end{equation*}\bigskip

\subsection{The Jacobi operator of minimal surfaces}
In this section, we recall some classical definitions and facts
about the Jacobi operator of minimal surfaces. Let  $M^2
\looparrowright \Mh^3$ be an orientable minimal surface immersed
into an oriented Riemannian manifold $(\Mh,\gh)$. Let $N_M$ be a
unit normal field along $M$, $A_M$  the second fundamental form of
the immersion with respect to the normal $N_M$, and  let $\rh$ be
the Ricci curvature of $\Mh$. The second variation of the volume
functional gives rise to the \emph{Jacobi operator} $J_M$ of $M$
(see \cite{L80})
\begin{equation}\label{E-jac0}
J_M :=-\De_M - (|A_M|^2 + \rh(N_M)),
\end{equation}
where $\De_M$ is the non-positive Laplacian on $M$ for the induced
metric.\medskip

Given a relatively compact regular domain $\Om$ on the surface $M$,
we let $\mathrm{Ind}(\Om)$ denote the number of negative eigenvalues
of $J_M$ for the Dirichlet problem in $\Om$. The \emph{Morse index}
of $M$ is defined to be the supremum
$$
\mathrm{Ind}(M):=\sup \{ \mathrm{Ind}(\Om); \Om \Subset M \}\le
\infty,
$$
taken over all relatively compact regular domains. Let $\la_1(\Om)$
be the least eigenvalue of the operator $J_M$ with the Dirichlet
boundary conditions in $\Om$. We call a relatively compact regular
domain $\Om$ \emph{stable} if $\la_1(\Om)>0$, \emph{unstable} if
$\la_1(\Om)<0$, and \emph{stable-unstable} if $\la_1(\Om)=0$. More
generally, we say that a domain $\Om$ (not necessarily relatively
compact) is \emph{r-stable} if any relatively compact subdomain is
stable. In the following proposition, we collect classical results
which will be used later on.

\begin{prop} \label{prop}
Given a minimal immersion $M^2 \looparrowright \Mh^3$, the following
properties hold.
\begin{enumerate}
\item[(i)] Let $\Om$ be a stable-unstable relatively compact domain. Then,
any smaller domain is stable while any larger domain is unstable.

\item[(ii)] We refer to the solutions of the equation $J_M(u)=0$ as
\emph{Jacobi functions} on $M$. Let $X_a: M^2\looparrowright
(\Mh^3,\gh)$ be a one-parameter family of oriented minimal
immersions, with variation field $V_a = \frac{\partial X_a}{\partial
a} $ and with unit normal $N_a$. Then, the function $\gh(V_a, N_a)$
is a Jacobi function on $M$.

\item[(iii)] Let $\Om$ be a relatively compact domain on a minimal submanifold $M$.
If there exists a positive function $u$ on $\Om$ such that
$J_M(u)\ge 0$, then $\Om$ is stable or stable-unstable.
\end{enumerate}
\end{prop}

\pf \emph{Assertion (i)} follows from the min-max characterization of
eigenvalues and the maximum principle. \emph{Assertion (ii)} appears in
\cite{BGS87} (Theorem 2.7 and its proof) in a more general
framework. For \emph{Assertion (iii)}, see the proof of Theorem~1 in
\cite{FS80}. \qed

%%%%%%%%%%%%%%%%%%%%%%%%%%%%%%%%%%%%%%%%%%%%%%%%%%%%%%%%

\section{Stable domains of revolution on the catenoids}\label{s-sd}

We consider a catenoid $\mathcal{C}$ given by the map,
$$
\mathcal{F}:\real\times[0,2\pi]\to\mathcal{C}\looparrowright \nil,
$$
$$
\mathcal{F}(t,\te)=(f(t)\cos \te, f(t)\sin\te, t),
$$
where $f$ is a global solution of (\ref{E-cat2}). It follows from
(\ref{tangents}) that the first fundamental form induced by
$\mathcal{F}$ is given by
\begin{equation*}\label{E-fff}
g_{\cF} =
\begin{pmatrix}
  1+f_t^2 & -\frac{1}{2}f^2 \\[5pt]
  -\frac{1}{2}f^2 & f^2 (1+\frac{1}{4}f^2)\\
\end{pmatrix}.
\end{equation*}
For later purposes, we introduce the functions
\begin{equation}\label{E-fff1}
G= f^2 (1+\frac{1}{4}f^2) \text{~~and~~} D =
\sqrt{\mathrm{Det}(g_{\cF})} = f \big( 1 + f_t^2 + \frac{1}{4} f^2
f_t^2 \big)^{1/2}.
\end{equation}

Let $N$ be a unit normal field to $\cF$. Writing $N=\al X + \be Y
+\ga Z$, we can choose $N$ to be
\begin{equation}\label{E-nw}
\left\{%
\begin{array}{lll}
\al & = & W ( -\cos \te -\frac{1}{2}f f_t \sin\te ),\\[4pt]
\be & =&  W ( -\sin\te + \frac{1}{2}f f_t \cos \te ),\\[4pt]
\ga &= & W f_t, \text{~~where}\\[4pt]
W &= & \big( 1+ f_t^2 + \frac{1}{4} f^2 f_t^2 \big)^{-1/2}.
\end{array}%
\right.
\end{equation}

\subsection{Jacobi functions coming from ambient Killing fields.}

Since the set $\{\xi, \eta, \zeta, \rho \}$ is a basis of Killing vector
fields, it follows from Proposition~\ref{prop}(ii) that the functions
\begin{equation}\label{E-nwj}
\left\{%
\begin{array}{lllll}
v_\xi&=&\gh(\xi,N)&=& W ( -\cos \te +\frac{1}{2}f  f_t \sin\te ),\\[4pt]
v_\eta&=&\gh(\eta,N)&=& W ( -\sin\te -\frac{1}{2}f f_t \cos \te ),\\[4pt]
v_\zeta&=&\gh(\zeta,N)&=& W f_t,\\
\end{array}%
\right.
\end{equation}
are  Jacobi functions on the surface $\cF$ (note that
$v_\rho=\gh(\rho,N)= 0$). \bigskip

\textbf{Remark}. The Jacobi functions $v_{\xi}, v_{\eta}$ and
$v_{\zeta}$ are linearly independent.

\subsection{A Jacobi function coming from the variation of the
family}\label{ss-vjf}

We now consider the one-parameter family of catenoids
$\{\mathcal{C}_a, a> 0 \}$, associated with the family of maps
\begin{equation}
\mathcal{F}(a,t,\te) = \big( f(a,t)\cos\te, f(a,t)\sin \te, t \big),
\quad a>0,
\end{equation}
where $f(a,\cdot)$ is the unique global solution of the Cauchy
problem (\ref{cp}). The variational field of this family is given by
\begin{equation}
\mathcal{F}_a(a,t,\te) = f_a(a,t)\cos\te \, X + f_a(a,t)\sin \te \,
Y,
\end{equation}
where $f_a(a,t):=\frac{\partial f}{\partial a}(a,t)$. By
Proposition~\ref{prop}(ii), this yields another Jacobi function on
$\mathcal{C}_a$, namely, $e(a,\cdot)= -\gh(\mathcal{F}_a,N).$ More
precisely,
\begin{equation}\label{e-jacobi}
e(a,t)= \big( W f_a \big)(a,t),
\end{equation}
where the function $W$ is given by the last line in (\ref{E-nw}). We
note that $e(a,\cdot)$ does not depend on $\te$ and is an even
function of $t$. Furthermore, since $f(a,0)=a$ and $f_t(a,0)=0$,
$\forall a>0$, we have $e(a,0)=1$, $\forall a>0.$
\bigskip

The rotationally invariant stable domains of the catenoids $\cC_a$
are described in the following theorem.

\begin{theorem}\label{t1}
Let $\mathcal{C}_a$ be a catenoid in $\nil$.  Then
\begin{enumerate}
\item[(i)] The upper (resp. the lower) half catenoid $\mathcal{C}_{a,+} =
\mathcal{C}_a\cap \{z>0\}$ (resp. $\mathcal{C}_{a,-} = \mathcal{C}_a
\cap \{z<0\}$) is r-stable.
\item[(ii)] The function $e(a,\cdot)$ is even and has exactly one zero $z(a)$ on $(0,\infty)$.
The domain $\mathcal{F}(a, [-z(a), z(a)], [0,2\pi])$ is a
stable-unstable domain in $\mathcal{C}_a$.
\item[(iii)] Given any $t_1 > 0$, there exists some $t_2 >0 $ such that the domain
$\mathcal{D}_a(-t_1,t_2)= \cF(a,[-t_1,t_2],[0,2\pi])$ is
stable-unstable. This implies in particular that both
$\mathcal{C}_{a,+}$ and $\mathcal{C}_{a,-}$ are maximal r-stable
rotationally invariant domains (\ie in some sense, stable-unstable).
\end{enumerate}
\end{theorem}

\pf \emph{Assertion (i)}. It follows from Section~\ref{ss-qa} that the
Jacobi function $v_{\zeta}$ is positive on $(0, + \infty)$ and
negative on $(- \infty , 0)$. The assertion follows from
Proposition~\ref{prop}(iii). \smallskip

\emph{Assertion (ii)}. We already know that $e(a,\cdot)$ is an even
function of $t$ and that $e(a,0)=1$ for all $a> 0$. \emph{Claim 1}.
The function $e(a,\cdot)$ has at most one zero in $(0, + \infty)$.
If not, $e(a,\cdot)$ would have two consecutive positive zeroes, $0
< z_1(a) < z_2(a)$ and the domain $\cF(a,[z_1(a), z_2(a)],[0,
2\pi])$ would be stable-unstable. According to
Proposition~\ref{prop}(i), this would contradict the $r$-stability of
$\mathcal{C}_{a,+}$ in Assertion~(i). \emph{Claim 2}. The function
$e(a,\cdot)$ has at least one zero in $(0, + \infty)$. Indeed,
$e(a,\cdot)$ has the sign of $f_a(a,t)$. Using the function $\phi$
defined by (\ref{E-phi}), we have
$$
\phi(a,f(a,t))\equiv t \text{~and~} \phi_a \big( a, f(a,t) \big) +
f_a(a,t) \, \phi_{\tau} \big(a, f(a,t) \big) \equiv 0
$$
for all $a, t > 0$. Since $\phi_{\tau}$ is positive, it suffices to
look at the sign of $\phi_a$. We find that
\begin{equation}\label{phi_a}
\phi_a(a,\tau) = \int_1^{\tau/a} \frac{a^2v^2 +
2}{\sqrt{(a^2v^2+4)(v^2-1)}} \, \dd v -
\frac{\tau}{2}\sqrt{\frac{\tau^2+4}{\tau^2-a^2}} %\leqno{(*)}
\end{equation}
and we easily conclude that $\phi_a(a,\tau)$ is positive when $\tau$
is large enough. It follows that $e(a,t)$ is negative for $t$ large
enough so that it must vanish at least once in $(0, + \infty)$.
\smallskip

\emph{Assertion (iii)}. Fix some $t_1 > 0$ and consider the function
$$
w(a,t_1,t) = v(a,t_1)e(a,t) + e(a,t_1) v(a,t),
$$
where we have written $v(a,t)$ instead of $v_{\zeta}(a,t)$ for
short. This is a Jacobi function on $\mathcal{C}_a$, which vanishes
at $t=-t_1$. Note that $w(a,t_1,0) = v(a,t_1) > 0$ because
$e(a,0)=1$ and $v(a,t) > 0$ for any $t>0$. As in the proof of
Assertion~(ii), Claim 1, we see that $w(a,t_1,\cdot)$ can vanish at
most once in $(- \infty,0)$ and $(0,\infty)$. It follows that
$w(a,t_1,\cdot)$ has exactly one zero in $(- \infty,0)$ -- namely
$-t_1$ -- and that it vanishes in $(0,\infty)$ if and only if it
takes some negative value
near infinity. %\smallskip
Recall that
\begin{equation}\label{E-v}
v(a,t) = \frac{f_t}{\sqrt{1 + f_t^2 + \frac{1}{4} f^2 f_t^2}}(a,t).
%\leqno{(a)}
\end{equation}

As in the proof of Assertion (ii), Claim 2, we use the functional
equations $\phi \big( a,f(a,t) \big) \equiv t$ and $\phi_{\tau}
\big( a,f(a,t) \big) \, f_t(a,t) \equiv 1$ for all $t>0$. Plugging
these relations into (\ref{E-v}), we find that
$$
v(a,t) = \tilde{v}\big( a, f(a,t) \big), ~~\forall t > 0,
%\leqno{(b)}
$$
where $\tilde{v}(a,\tau) = \big( 1 + \frac{\tau^2}{4} +
\phi_{\tau}^2(a,\tau)\big)^{-1/2}$. Similar computations yield the
relation
$$
e(a,t) = \tilde{e}\big( a, f(a,t) \big), ~~~\forall t > 0,
%\leqno{(c)}
$$
where $\tilde{e}(a,\tau) = - \phi_a(a,\tau) \tilde{v}(a,\tau)$.
Define $\tau_1 := f(a,t_1)$ and
$$
\tilde{w}(a,\tau_1,\tau) = \tilde{v}(a,\tau_1) \tilde{e}(a,\tau) +
\tilde{e}(a,\tau_1) \tilde{v}(a,\tau), %\leqno{(d)}
$$
so that $w(a,t_1,t) = \tilde{w}\big(a,\tau_1,f(a,t)\big)$. Then,
$$
\tilde{w}(a,\tau_1,\tau) =  -\tilde{v}(a,\tau)
\tilde{v}(a,\tau_1)\big( \phi_a(a,\tau) + \phi_a(a, \tau_1) \big).
%\leqno{(e)}
$$
Using (\ref{phi_a}), we conclude that $w$ is negative when $\tau$ approches
infinity, for any given $a, t_1 > 0$. This proves the existence of a
positive $t_2$ such that the domain $\mathcal{D}_a(-t_1,t_2)$ ist
stable-unstable. The last assertion follows immediately. \qed
\bigskip

\textbf{Remarks}. (i) Consider the family of curves $\Gamma_a : t
\mapsto \big( f(a,t),t \big)$. This family admits an envelope $\cE$ and
the values $\pm z(a)$ correspond to the points at which the curve
$\Gamma_a$ is tangent to $\cE$. (ii) Using (\ref{E-nw}) and
Section~\ref{ss-qa}, we can
see that the Gauss map of the catenoid $\cC_a$ covers a closed
symmetric strip about the equator of the unit sphere in the Lie
algebra $\mathcal{L}(\nil)$. This strip, whose width depends on $a$,
is strictly contained in the sphere minus the south and north poles.
Each point of the open strip is covered exactly twice, except the
points of the equator which are covered once (look at the variations
of the $Z$-component $\gamma$ of the vector $N$).

%%PRIVATE

%%%%%%%%%%%%%%%%%%%%%%%%%%%%%%%%%%%%%%%%%%%%%%%%%%%%%%%%%%%%%%%%%%

\section{The index of the catenoids $\mathcal{C}_a$ in $\nil$}\label{s-ind}

In this section, we study the Morse index of the catenoids $\cC_a$.
It turns out that the representation $\cF$ given by (\ref{param}),
with the function $f$ satisfying (\ref{E-cat2}), is not well-adapted
to Fourier analysis on $\mathcal{C}_a$ because the vectors $\cF_t$
and $\cF_{\theta}$ are not orthogonal. To avoid this problem, we
introduce a perturbed representation,
\begin{equation*}\label{pparam}
\widetilde{\mathcal{F}}(t,\te) :=\mathcal{F}(t,\te+\varphi(t))=
\Big(f(t)\cos(\te +\varphi (t)), f(t)\sin (\te +\varphi (t)),
t\Big).
\end{equation*}
The tangent vectors are given by
\begin{equation*}
\left\{%
\begin{array}{lcl}
\widetilde{\mathcal{F}}_t(t,\te)&=& \mathcal{F}_t(t,\te +
\varphi(t))+ \varphi_t (t)\mathcal{F}_\te(t,\te+\varphi(t)),\\[4pt]
\widetilde{\mathcal{F}}_\te(t,\te)&=&  \mathcal{F}_{\te}(t,\te + \varphi(t)).\\
\end{array}%
\right.
\end{equation*}

It follows that the representation $\widetilde{\mathcal{F}}$ is
orthogonal -- \ie the vectors $\widetilde{\mathcal{F}}_t$ and
$\widetilde{\mathcal{F}}_{\theta}$ are orthogonal -- if and only if
the function $\varphi$ satisfies the differential equation
\begin{equation}\label{varphi}
\varphi_t =\frac{2}{4+f^2}.
\end{equation}

From now on, we choose $\varphi$ to be the solution of
(\ref{varphi}) such that $\varphi(0)=0$. \smallskip

Note that in the above expressions, we have omitted the dependence
on the parameter $a$. The unit normal vector to $\mathcal{C}_a$ at
the point $\widetilde{\mathcal{F}}(t,\theta)$ is
$\widetilde{N}(t,\te) = N(t,\te + \varphi(t))$. In the
representation $\widetilde{\mathcal{F}}$, the Riemannian metric
induced by the immersion $\cC_a \looparrowright \nil$ is of the form
$D^2G^{-1} \dd t^2 + G \dd \theta^2$, with the functions $D, G$ as
in (\ref{E-fff1}). It follows that the Laplacian on $\mathcal{C}_a$
is given, in the representation $\widetilde{\mathcal{F}}$, by the
expression
\begin{equation*}\label{laplacian}
\widetilde{\Delta} = \frac{1}{D}\,
\partial_t \Big( \frac{G}{D} \partial_t
\Big) + \frac{1}{G} \, \partial^2_{\theta \theta}.
\end{equation*}
We introduce the operator
\begin{equation*}\label{E-L}
\widetilde{L} = - \frac{1}{D}\, \partial_t \Big( \frac{G}{D}
\partial_t \Big),
\end{equation*}
and the function
\begin{equation*}\label{E-V}
\widetilde{V} = \big( \rich(\widetilde{N}) + |\widetilde{A}|^2\big),
\end{equation*}
which only depend on the variable $t$ (and the parameter $a$). In
the parametrization $\widetilde{\mathcal{F}}$, the Jacobi operator
(\ref{E-jac0}) of the immersion $\cC_a \looparrowright \nil$ is
given by the expression
\begin{equation*}\label{E-Jact}
\widetilde{J} = \widetilde{L} - \widetilde{V} - \frac{1}{G} \,
\partial^2_{\theta \theta}.
\end{equation*}

%%PRIVATE

We have the following lemma.

\begin{lem}\label{L-t2-1}
With the above notations, the function  $\widetilde{V}$ on the
catenoid $\cC_a$ is given by,
\begin{equation*}\label{E-V1}
\widetilde{V} = \frac{2a^2}{f^4} + \frac{2(a^2+4)}{(4+f^2)^2}.
\end{equation*}
Furthermore, the function $G\widetilde{V}$ is equal to
$\frac{a^2}{2} \, \frac{4+f^2}{f^2} + \frac{a^2+4}{2} \,
\frac{f^2}{4+f^2}$ and satisfies the inequalities
\begin{equation*}\label{E-V2}
(a^2+2) \sqrt{1 - \frac{4}{(a^2+2)^2}} = a \sqrt{a^2+4}\le
(G\widetilde{V})(a,t) \le a^2+2,
\end{equation*}
for all $a>0$ and all $t \in \R$.
\end{lem}

\pf For the catenoid $\cC_a$, the function $f$ satisfies the
differential equations (\ref{1st}) and (\ref{E-cat2}) and we have $W
= \frac{a}{f}$, where the function $W$ is defined in (\ref{E-nw}).
The $Z$-component $\gamma$ of the unit normal $\widetilde{N}$ is a
Jacobi function which only depends on $t$, hence
$\widetilde{L}(\gamma) = \widetilde{V}\gamma$. Using (\ref{1st}) and
(\ref{E-cat2}) again, we can compute $\widetilde{L}(\gamma)$ and
derive the formulas for $\widetilde{V}$ on the catenoid $\cC_a$. The
second assertion follows easily. \qed \bigskip

%%PRIVATE

Let $\widetilde{v}_{\xi}$ and $\widetilde{v}_{\eta}$ be the
expressions of the Jacobi functions associated with the Killing
fields $\xi$ and $\eta$ in the parametrization $\widetilde{\cF}$. It
follows from (\ref{E-nwj}) that
$$
\widetilde{v}_{\xi}(t,\theta) = \gh \Big(
\xi(\widetilde{\mathcal{F}}(t, \te)), \widetilde{N}(t,\theta) \Big)
= W \Big( -\cos (\te + \varphi) +\frac{1}{2}f f_t \sin (\te +
\varphi) \Big),
$$
and similarly for $\widetilde{v}_{\eta}$ (we have omitted the
dependence on $a$). We introduce the smooth function $\psi (a,t)$
such that
\begin{equation*}\label{E-psi}
\left\{%
\begin{array}{rcl}
\cos \psi &=& (1 + \frac{1}{4}f^2 f_t^2)^{-1/2}, \\[4pt]
\sin \psi &=& \frac{1}{2}f f_t(1 + \frac{1}{4}f^2 f_t^2)^{-1/2},\\[4pt]
\psi(a,0) &=& 0.
\end{array}%
\right.
\end{equation*}

It follows immediately that
\begin{equation*}\label{E-jac2}
\left\{%
\begin{array}{rcl}
\widetilde{v}_{\xi}(a,t,\theta) & = & - W_1(a,t) \cos \big( \theta +
\varphi(a,t)
+ \psi(a,t) \big), \\[4pt]
\widetilde{v}_{\eta}(a,t,\theta) & = & - W_1(a,t) \sin \big( \theta +
\varphi(a,t) + \psi(a,t) \big), \text{ ~where}\\[4pt]
W_1 & = & W (1 + \frac{1}{4}f^2 f_t^2)^{1/2} = \frac{1}{f} \,
\sqrt{\frac{4a^2+f^4}{4+f^2}}.
\end{array}%
\right.
\end{equation*}

With the above notations, we have the following lemma.

\begin{lem}\label{L-t2-2}
Let $\omega := \varphi + \psi$, a function of the variable $t$ and
the parameter $a$. Then,
\begin{enumerate}
    \item[(i)] The functions
    \begin{equation}\label{E-jac4}
    \left\{%
    \begin{array}{ll}
    w_1(a,t,\theta ) := W_1(a,t) \cos (\omega(a,t)) \cos \theta,  \\
    w_2(a,t,\theta ) := W_1(a,t) \cos (\omega(a,t)) \sin \theta,  \\
    w_3(a,t,\theta ) := W_1(a,t) \sin (\omega(a,t)) \cos \theta, \\
    w_4(a,t,\theta ) := W_1(a,t) \sin (\omega(a,t)) \sin \theta, \\
    \end{array}%
    \right.
    \end{equation}
    are bounded Jacobi functions on $\cC_a$, $\widetilde{J}(w_i)=0$, for $1 \le i \le 4$.
    \item[(ii)] The function $\omega(a,\cdot)$ is an odd function
    of $t$, satisfying $\omega(a,0)=0$ and $\omega_t = 4 f^2 (f^4 + 4 a^2)^{-1}$.
    \item[(iii)] Let $\Omega(a) := \lim_{t \to + \infty} \omega(a,t)$. Then
    \begin{equation*}\label{E-crit1}
    \Omega(a) = 2a \int_{a}^{\infty} \frac{u^2 \sqrt{u^2+ 4}}{(u^4+4a^2)
    \sqrt{u^2-a^2}} \, \dd u.
    \end{equation*}
    \item[(iv)] For all $a>0$, we have $\frac{\pi}{2} < \Omega(a) \le \pi$ and
    the lower and upper bounds are achieved as limits when $a$ tends respectively
    to zero and infinity.
\end{enumerate}
\end{lem}

\pf  \emph{Assertion (i)} follows from the equalities
$\widetilde{v}_{\xi} = - w_1 + w_4$, $\widetilde{v}_{\eta} = - w_2 -
w_3$, and the fact that the operator $\widetilde{J}$ separates
variables. \emph{Assertion (ii)}. The computation of $\omega_t$ is
straightforward. To prove \emph{Assertion (iii)}, we use the fact that
$f_t$ is positive for positive $t$ and can be computed from
(\ref{1st}), namely,
$$
f_t = \frac{2 \sqrt{f^2-a^2}}{a \sqrt{f^2+4}}.
$$
We write
$$
\omega_t = \frac{2a f^2 \sqrt{f^2+4}}{(f^4 + 4 a^2) \sqrt{f^2-a^2}}\, f_t
$$
for $t > 0$, and we compute the integral $\int_0^t \omega_{\tau} \,
\dd \tau$ by making the change of variables $u = f(t)$.
\emph{Assertion (iv)}. Assume by contradiction that  $\Omega(a_0) >
\pi$ for some $a_0$. There would then exist a value $t_0$ such that
$\omega(a_0,t_0) = \pi$. The function $w_3$ in (\ref{E-jac4}) would
then vanish on the circles $\widetilde{\cF}(a_0,\{0\},[0, 2\pi])$
and $\widetilde{\cF}(a_0,\{t_0\},[0, 2\pi])$. Because this function
is a Jacobi function, this would contradict Assertion (i) in
Theorem~\ref{t1}. The fact that $\frac{\pi}{2} < \Omega(a)$ follows
from a direct estimate of the integral, \cite{Li10}. Indeed, making
the change of variables $u = av$, we get $\Omega (a) =
\Omega_1(4/a^2)$ where
$$\Omega_1(b) = 2 \int_1^{\infty} \frac{v^2 \sqrt{v^2+b}}{(v^4+b)\sqrt{v^2-1}} \,
\dd v = \int_1^{\infty} \frac{\sqrt{u}}{\sqrt{u-1}} \,
\frac{\sqrt{u+b}}{u^2+b} \, \dd u > I(b),$$
where
$$I(b) = \int_1^{\infty} \frac{\sqrt{u+b}}{u^2+b} \, \dd u$$
and we claim that $I(b) > \frac{\pi}{2}$. To prove this last
assertion, we consider two cases, $0 \le b \le 1$ and $b > 1$.
\medskip

\noindent $\bullet~$ We have $I(0)=2$ and for $0 < b \le 1$,
$$
I(b) > \int_1^{\infty} \frac{\sqrt{u}}{u^2+1} \, \dd u \ge 2
\int_1^{\infty} \frac{1}{v^2+1} \, \dd v = \frac{\pi}{2}.
$$

\noindent $\bullet~$ When $b \ge 1$, we can write
$$I(b) = \int_1^b  \frac{\sqrt{u+b}}{u^2+b} \, \dd u + \int_b^{\infty}
\frac{\sqrt{u+b}}{u^2+b} \, \dd u$$
and estimate the integrals on the right-hand side separately.
$$
\begin{array}{lcl}
\displaystyle \int_b^{\infty} \frac{\sqrt{u+b}}{u^2+b} \, \dd u & >&
\displaystyle \int_b^{\infty} \frac{\sqrt{u}}{u^2+b} \, \dd u = 2
\displaystyle \int_{\sqrt{b}}^{\infty} \frac{v^2}{v^4+b} \, \dd v \\[6pt]
& \ge& 2 \displaystyle \int_{\sqrt{b}}^{\infty} \frac{1}{v^2+1} \,
\dd v = \pi - 2 \arctan \sqrt{b}.\\
\end{array}
$$

On the other hand,
$$
\begin{array}{lcl}
\displaystyle \int_1^b  \frac{\sqrt{u+b}}{u^2+b} \, \dd u & \ge&
\sqrt{b} \displaystyle \int_1^b  \frac{1}{u^2+b} \, \dd u =
\displaystyle \int_{\frac{1}{\sqrt{b}}}^{\sqrt{b}} \frac{dv}{v^2+1}\\[8pt]
&  =& \arctan (\sqrt{b}) - \arctan (\frac{1}{\sqrt{b}})\\[6pt]
&  =& 2 \arctan
(\sqrt{b}) - \frac{\pi}{2}.\\
\end{array}
$$
It follows that $I(b) > \frac{\pi}{2}$ and hence that $\Omega(a) >
\frac{\pi}{2}$. Recall that $\Omega(a) \le \pi$ for geometric
reasons. Clearly, when $b$ tends to zero, $\Omega_1(b)$ tends to
$\pi$, and hence $\Omega(a)$ tends to $\pi$ when $a$ tends to
infinity. Making the change of variable $u = \sqrt{b} \, v$, one can
show that $I(b)$ tends to $\frac{\pi}{2}$ when $b$ tends to
infinity. On the other hand, it is easy to see that $\Omega_1(b) -
I(b)$ tends to $0$ when $b$ tends to infinity. It follows that
$\Omega(a)$ tends to $\frac{\pi}{2}$ when $a$ tends to zero. This
finishes the proof of the lemma. \qed
\bigskip

%%PRIVATE

\begin{lem}\label{L-t2-3}
For $k \in \N$, consider the operator $\widetilde{L_k} :=
\widetilde{L} + \frac{k^2}{G} - \widetilde{V}$ in $L^2([-r,r], D\,
\dd t)$, with Dirichlet boundary conditions. Then,
\begin{enumerate}
    \item[(i)] For any $r>0$, the operator $\widetilde{L_k}$ has at most one
    negative eigenvalue (with multiplicity one).
    \item[(ii)] For all $k \ge \sqrt{a^2+2}$ and $r>0$, the operator $\widetilde{L_k}$ is
    positive.
\end{enumerate}
\end{lem}\bigskip

\pf \emph{Assertion (i)}. Recall that the eigenvalues of a
Sturm-Liouville problem with Dirichlet boundary conditions are
always simple. If $\widetilde{L_k}$ had at least two negative
eigenvalues, we would have an eigenfunction $v$ of $\widetilde{L_k}$
associated with a negative eigenvalue and having one zero in
$(-r,r)$. The function $v \cos(k\theta)$ would be an eigenfunction
of the Jacobi operator $\widetilde{J}$ with negative eigenvalue,
vanishing on the boundary of an annulus contained in $\cC_{a,+}$ or
in $\cC_{a,-}$. This would contradict Assertion (i) in
Theorem~\ref{t1}. \emph{Assertion (ii)}. By Lemma~\ref{L-t2-1}, $G
\widetilde{V} \le a^2+2$ and the second assertion follows from the
positivity of the operator $\widetilde{L}$ in $L^2([-r,r], D\, \dd
t)$. \qed

\newpage

\begin{thm}\label{t2}
Consider the catenoids $\{\cC_a, a > 0\}$ in $\nil$.
\begin{enumerate}
    \item[(i)] For all $a>0$, the catenoid $\cC_a$ has finite Morse index at least
    equal to $3$ and at most equal to $1 + 2[\sqrt{a^2+2}]$, where $[x]$ is the
    integer part of $x$. In particular, the index of $\cC_a$ is equal to $3$ for
    $a$ close to zero.
    \item[(ii)] When $a$ tends to infinity, the index of $\cC_a$ grows
    at least like $\sqrt{3}\, a$. In particular, it tends to infinity
    when $a$ tends to infinity.
\end{enumerate}
\end{thm}

\pf Fourier analysis and Lemma~\ref{L-t2-3}(i) show that the Morse
index of $\cC_a$ is equal to $1$ plus twice the number of positive
integers $k$ such that the operator $\widetilde{L}_k$ has a negative
eigenvalue. \emph{Assertion (i)}. The fact that the index of $\cC_a$
is at most $1 + 2[\sqrt{a^2+2}]$ follows from Lemma~\ref{L-t2-3}(ii).
By Lemma~\ref{L-t2-2}(iv), for any $a>0$, $\Omega(a) > \pi/2$. Since
$\omega(0)=0$, there exists some $r_a > 0$ such that
$\omega(a,r_a)=\frac{\pi}{2}$. The functions $w_1, w_2$ of
Lemma~\ref{L-t2-2} (i) are Jacobi functions which vanish on the
boundary of the domain $\cF(a,(-r_a,r_a),[0,2\pi])$. It follows
easily that the index of the operator $\widetilde{L}_1$ is equal to
$1$ and hence the index of the catenoid $\cC_a$ is at least $3$.
\emph{Assertion (ii)}. To determine whether the index of
$\widetilde{L}_k$ is $1$ or $0$, consider the associated quadratic
form on functions $\phi \in C_0^1(\R)$,
$$
Q_k(\phi) = \int_{- \infty}^{\infty} \big\lbrace \frac{G}{D}
\phi_t^2 + (k^2 - G \widetilde{V}) \frac{D}{G} \phi^2 \big\rbrace \,
dt.
$$
Write $\phi(t) = \psi\big(s(t)\big)$, with
$$
s_t = \frac{D}{G} = \frac{4}{a(4+f^2)}, ~ s(0) = 0.
$$
The function $s$ is a diffeomorphism from $\R$ onto $\big(-S(a),
S(a) \big)$, where
\begin{equation}\label{E-ig1}
S(a) = 2 \int_0^{\infty} \frac{f_t \, dt}{\sqrt{(4+f^2)(f^2-a^2)}} =
\frac{2}{a} \int_1^{\infty} \frac{du}{\sqrt{(u^2 +
\frac{4}{a^2})(u^2-1)}}.
\end{equation}
It follows that
$$
Q_k(\phi) = \int_{-S(a)}^{S(a)} \big\lbrace \psi_s^2 + \big(k^2 -
U(s) \big) \psi^2 \big\rbrace \, ds,
$$
where the function $U$ is defined by $U\big( s(t) \big) = (G
\widetilde{V})(t)$. Choose the function $\psi$ to be $\psi_0(s) =
\cos \big( \frac{\pi s}{2S(a)} \big)$ and let $\phi_0$ be the
corresponding function. Using Lemma~\ref{L-t2-1}, one finds that
$Q_k(\phi_0) < 0$, \ie that the index of $\widetilde{L}_k$ is $1$,
as soon as
\begin{equation}\label{E-ig2}
k^2 < (a^2+2)\sqrt{1 - \big(\frac{2}{a^2+2}\big)^2} - \big(
\frac{\pi}{2 S(a)} \big)^2 .
\end{equation}
By (\ref{E-ig1}), $S(a) = \frac{\pi}{a} - \frac{4}{a^3} J(a)$, where
the function $J(a)$ is given by
$$
J(a) = 2 \int_1^{\infty} \frac{dv}{v (v + \sqrt{v^2 +
\frac{4}{a^2}}) \sqrt{(v^2-1)(v^2 + \frac{4}{a^2})}}.
$$
This function tends to $\frac{\pi}{4}$ when $a$ tends to infinity
and hence the right-hand side of (\ref{E-ig2}) is equivalent to
$\frac{3a^2}{4}$ when $a$ tends to infinity. This proves the second
assertion. \qed \medskip

\newpage

\textbf{Remarks}.
\begin{enumerate}
    \item[(i)] Given $a>0$, there is a simple criterion to decide whether
    the operator $\widetilde{L}_k$ has a
    negative eigenvalue in the interval $[-r,r]$ (with Dirichlet
    boundary conditions). Let $u_k$ be the solution of the Cauchy problem
    $\widetilde{L}_k(u)=0$, $u(0)=1$ and $u_t(0)=0$. If $u_k$ has a
    zero in the interval $(0,r)$, then $\widetilde{L}_k$ has a
    negative eigenvalue in $[-r,r]$; if $u_k$ does not vanish in the
    interval $(0,r)$, then $\widetilde{L}_k(u) \ge 0$ in $[-r,r]$.
    \item[(ii)] Using the fact that the metric $\gh$ on $\nil$ is
    left-invariant, one can easily express the associated Levi-Civita
    connexion and curvature tensors in the orthonormal basis
    $\{X,Y,Z\}$ of left-invariant vector fields. In particular, given a
    unit vector $N = \alpha X + \beta Y + \gamma Z$, we find the
    following formula for the Ricci curvature,%
    $$\rich(N,N) = - \frac{1}{2} + \gamma^2.$$
    %
%%PRIVATE
    %
    \item[(iii)] Using the preceding remark, we can write the Jacobi
    operator on an orientable minimal surface in $\nil$ as
    $$J = - \Delta + \frac{1}{2} - \gamma^2 - |A|^2,$$
    where $\gamma$ is the $Z$-component of the unit normal to the
    surface. Using the fact that the scalar curvature of $\nil$ is
    $-\frac{1}{4}$, we also have the formula
    $$J = - \Delta + \frac{1}{4} + K_M - \frac{1}{2}|A|^2,$$
    where $K_M$ is the Gauss curvature of the surface $M$.
    \item[(iv)] Using Lemma~\ref{L-t2-1} and the second remark, we
    deduce the following expression for the second fundamental form
    of the catenoid $\cC_a$ in $\nil$,
    $$|A|^2 = \frac{1}{2} - \frac{4}{f^2} + \frac{4(a^2+4)}{f^2(f^2+4)}
    + \frac{2(a^2+4)}{(f^2+4)^2}.$$
    This shows that the norm squared of the second fundamental form tends to
    $\frac{1}{2}$ uniformly at infinity. This is in contrast with the
    situation in $\R^3, \HH^2\times \R$ or $\HH^3$.
\end{enumerate}

%%%%%%%%%%%%%%%%%%%%%%%%%%%%%%%%%%%%%%%%%%%%%%%%%%%%%%%%%

\section{Catenoids in higher dimensions}\label{s-hd}

In this section, we study the rotationally symmetric stable domains
on the higher dimensional catenoids. Let $\niln$ be the
$(2n+1)$-dimensional Heisenberg group. As in Section \ref{s-pc}, we
use the exponential coordinates and choose the left-invariant metric
$\gh$ to be such that the left-invariant vector fields $\{X_1, \cdots ,
X_n, Y_1, \cdots, Y_n, Z\}$ form an orthonormal basis, where
\begin{equation*}\label{E-hd-1}
\left\{%
\begin{array}{rcl}
X_i(x,y,z) &=& \partial_{x_i} - \frac{1}{2}y_i \partial_z, ~ 1 \le i \le n,\\[4pt]
Y_i(x,y,z) &=& \partial_{y_i} + \frac{1}{2}x_i \partial_z, ~ 1 \le i \le n,\\[4pt]
Z(x,y,z)   &=& \partial_z.
\end{array}%
\right.
\end{equation*}

We look for hypersurfaces of revolution of the form
\begin{equation*}\label{E-hd-2}
\cF :
\left\{%
\begin{array}{l}
\R \times S^{2n-1} \to \niln,\\
(t,\theta) \mapsto \cF(t,\theta) = \big( f(t) \theta ,t \big),\\
\end{array}%
\right.
\end{equation*}
where $f$ is a positive function of $t$. If follows from \cite{Fi96,
FMP99} that the hypersurface $\cF$ is minimal if and only if $f$
satisfies the second order differential equation,
\begin{equation*}\label{E-hd-3}
f (4+f^2) f_{tt} = 4(2n-1) (1+f_t^2) + (2n-2) f^2 f_t^2.
\end{equation*}

As in Section~\ref{ss-qa}, one can show that for $a>0$, there is a
unique maximal solution $f(a,t)$ such that $f(a,0)=a$ and
$f_t(a,0)=0$. This is an even function of $t$ defined on the
interval $(-T(a),T(a))$, where $T(a)$ is finite when $n \ge 2$. As
in dimension $3$ ($n=1$), the above differential equation admits a
first integral,
\begin{equation*}\label{E-hd-4}
f^{2n-1} \, \big( 1 + f_t^2 + f^2 f_t^2 \big)^{-1/2} \equiv
a^{2n-1}.
\end{equation*}

As in (\ref{E-nw}), we let $W := \big( 1 + f_t^2 + f^2 f_t^2
\big)^{-1/2}$. We also use the following notations,
\begin{equation*}\label{E-hd-5}
\left\{%
\begin{array}{rcl}
\cC_a &=& \cF \Big(a, \big( -T(a),T(a) \big), S^{2n-1} \Big), \\[4pt]
\cC_{a,+} &=& \cF\Big( a, \big( 0,T(a) \big), S^{2n-1} \Big), \\[4pt]
\cC_{a,-} &=& \cF\Big( a, \big( -T(a), 0 \big), S^{2n-1} \Big), \\[4pt]
\cD_a(r,s) &=& \cF\Big( a, (r,s), S^{2n-1}\Big).
\end{array}%
\right.
\end{equation*}

We can now state the following result.

\begin{thm}\label{t3}
Assume that $n \ge 2$ and $a>0$.
\begin{enumerate}
    \item[(i)] The half-catenoids $\cC_{a,\pm}$ are r-stable.
    \item[(ii)] There exists some $z(a) > 0$ such that the domain
    $\cD_a(-z(a), z(a))$ is stable-unstable. In particular, the
    catenoid $\cC_a$ has index at least $1$.
    \item[(iii)] There exists some $\ell(a) > 0$ such that the domain
    $\cD_a(- \ell(a), T(a))$ is r-stable.
    \item[(iv)] For any $r > \ell(a)$, there exists some $s > 0$ such
    that the domain $\cD_a(-r,s)$ is stable-unstable.
\end{enumerate}
\end{thm}

\pf The proof relies on the expressions of two explicit Jacobi
functions on $\cC_a$, namely the Jacobi functions $v(a,t) =
\gh(N,Z)$, and $e(a,t) = - \gh(\cF_a,N)$, where $N$ is a unit normal
to $\cC_a$, and $\cF_a$ is the variation field along $\cF$ when the
parameter $a$ varies. As in dimension $2$, we have $v(a,t) = W(a,t)
f_t(a,t)$ and Assertion (i) follows immediately from the fact that
$f_t(a,t) > 0$ for $t > 0$. \medskip

To prove the other Assertions, notice that $e(a,t)$ is an even
function of $t$ which can be studied using the inverse function
$\phi(a,\tau)$ of the function $f(a,\cdot) : [0,\infty) \to [a,
T(a))$. It turns out that
\begin{equation*}\label{E-hd8}
\phi(a,\tau) = \frac{a^{2n-1}}{2} \int_a^{\tau}
\sqrt{\frac{u^2+4}{u^{4n-2} - a^{4n-2}}} \, \dd u.
\end{equation*}
This formula shows that $\phi(a,\tau)$ has a finite limit $T(a)$
when $\tau$ tends to infinity and that its derivative
$\phi_a(a,\tau)$ has a positive finite limit when $\tau$ tends to
infinity. We now use the same method as in the proof of
Theorem~\ref{t1}. Assertion (ii), follows from the fact that $e(a,0)=1$
and that $e(a,t)$ takes negative values near infinity. For the
proofs of Assertions~(iii) and (iv), we use the fact that in higher
dimensions ($n \ge 2$), both $\phi(a,\tau)$ and $\phi_a(a,\tau)$
have finite limits at infinity, so that the higher dimensional case
differs from the case in which $n=1$. \qed \medskip

\textbf{Remark}. Theorem~\ref{t1}(iii) tells us that the
half-catenoids $\cC_{a,\pm}$ in $\nil$ are stable-unstable, \ie that
they satisfy the \emph{Lindeloef's property} as defined in
\cite{BS08, BS10}. Theorem~\ref{t3}(iii) and (iv) tell us that
catenoids in $\niln$, $n \ge 2$, do not satisfy Lindeloef's
property. As for catenoids in $\R^{n+2}$ and $\HH^n\times\R$, $n\ge
2$, this is related to the fact that these catenoids have finite
height. \bigskip

%%PRIVATE

\bigskip
\bigskip

\begin{small}

\begin{tabular}{l}
Pierre B\'{e}rard\\
Universit\'{e} Grenoble 1\\
Institut Fourier (\textsc{ujf-cnrs})\\
B.P. 74\\
38402 Saint Martin d'H\`{e}res Cedex\\
France\\
\verb+Pierre.Berard@ujf-grenoble.fr+\\
\end{tabular}\bigskip

\begin{tabular}{l}
Marcos P. Cavalcante\\
Universidade Federal de Alagoas\\
Instituto de Matem\'{a}tica\\
57072-900 Macei\'o-AL\\
Brazil\\
\verb+marcos.petrucio@pq.cnpq.br+\\
\end{tabular}\

\end{small}
\end{document}